\documentclass[a4paper,fleqn]{cas-sc}
\usepackage[numbers]{natbib}
\usepackage{amsthm}
\usepackage{orcidlink}
\usepackage{algorithm}
\usepackage{algorithmic}

\newcommand{\algorithmicbreak}{\textbf{break}}
\newcommand{\BREAK}{\STATE{\algorithmicbreak}}
\newtheoremstyle{style}{}{}{\normalfont}{}{\bfseries}{.}{ }{}
\theoremstyle{style}
\newtheorem{theorem}{Theorem}[section]
\newtheorem{definition}{Definition}[section]
\newtheorem{lemma}{Lemma}[section]

\newtheorem{example}{Example}[section]
\def\tsc#1{\csdef{#1}{\textsc{\lowercase{#1}}\xspace}}
\tsc{WGM}
\tsc{QE}

\begin{document}
\let\WriteBookmarks\relax
\def\floatpagepagefraction{1}
\def\textpagefraction{.001}
\shorttitle{Relating bubble sort to birthday problem}
\shortauthors{Jichu Jiang}
\title [mode = title]{Relating bubble sort to birthday problem}
\author[1]{Jichu Jiang \orcidlink{0009-0008-9777-610X}}
\cormark[1]
\ead{saisyc@outlook.com}
\affiliation[1]{organization={Beihang University},
            city={Beijing},
            postcode={100191},
            country={China}}
\cortext[1]{Corresponding author}

\begin{abstract}
    Birthday problem is a well-known classic problem in probability theory widely applied in cryptography, and bubble sort is a popular sorting algorithm leading to some interesting theoretical problems in computer science. However, the relation between bubble sort and birthday problem has not been discovered. By relating bubble sort to birthday problem,  based on a generalization of Poisson limit theorem for dissociated random variables, this paper offers an intuitive explanation to naturally indicate that $\displaystyle \frac{n - P_{n}}{\sqrt{n}}$ converges in distribution to the standard Rayleigh distribution, where $P_{n}$ is the number of passes required to bubble sort $n$ distinct elements. Then asymptotic distributions and statistical characteristics of bubble sort and birthday problem are presented. Moreover, this paper discovers that some common optimizations of bubble sort could lead to average performance degradation.
\end{abstract}


\maketitle

\section{Introduction}

    Birthday problem is a well-known classic problem in probability theory \cite{The_matching_birthday_and_the_strong_birthday_problem:_a_contemporary_review}, invariably mentioned in both mathematics textbooks such as \cite{Probability:_Theory_and_Examples, Analytic_Combinatorics} and computer science textbooks such as \cite{Introduction_to_Algorithms_Third_Edition, The_art_of_computer_programming_volume_3:_(2nd_ed.)_sorting_and_searching}. Denote $C_{n}$ as the number of persons at the first collision in a year of $n$ days, then there holds $\displaystyle \mathbb{E} \left( C_{n} - 1 \right) = \sqrt{\frac{\pi n}{2}} \pm \Theta \left( 1 \right)$ \cite{Analytic_Combinatorics}, and the limiting distribution of $\displaystyle \frac{C_{n} - 1}{\sqrt{n}}$ is Rayleigh distribution \cite{Probability:_Theory_and_Examples}. Although birthday problem and its generalizations received much attention for their applications to cryptography such as \cite{A_nonuniform_birthday_problem_with_applications_to_discrete_logarithms, Asymptotic_results_for_the_number_of_Wagner’s_solutions_to_a_generalised_birthday_problem, Nonuniform_birthday_problem_revisited:_Refined_analysis_and_applications_to_discrete_logarithms}, sorting algorithms have not been related to birthday problem.
    
    Bubble sort is a popular sorting algorithm \cite{Bubble_sort:_an_archaeological_algorithmic_analysis}, leading to some interesting theoretical problems in computer science \cite{The_art_of_computer_programming_volume_3:_(2nd_ed.)_sorting_and_searching}. Denote $P_{n}$ as the number of passes required to sort $n$ elements, then there holds $\displaystyle \mathbb{E} \left( n - P_{n} \right) = \sqrt{\frac{\pi n}{2}} \pm \Theta \left( 1 \right)$ \cite{The_art_of_computer_programming_volume_3:_(2nd_ed.)_sorting_and_searching}, whereas the limiting distribution of $\displaystyle \frac{n - P_{n}}{\sqrt{n}}$ has not been discovered. Although researchers such as \cite{Faulttolerant_maximal_localconnectivity_on_Bubblesort_star_graphs, Conditional_diagnosability_of_bubblesort_star_graphs, Edge_faulttolerant_strong_Menger_edge_connectivity_of_bubblesort_star_graphs} focused on the structure of bubble sort graph in recent years, little attention has been paid to the asymptotic distribution of the number of the passes.
    
    Moreover, in computer engineering, researchers such as \cite{Enhancement_of_Bubble_and_Insertion_Sort_Algorithm_Using_Block_Partitioning, Smart_Bubble_Sort:_A_Novel_and_Dynamic_Variant_of_Bubble_Sort_Algorithm, Experiment_Analysis_on_the_Bubble_Sort_Algorithm_and_Its_Improved_Algorithms} have often compared the average performance of bubble sort with those of other sorting algorithms throughout many years, which is determined by the distribution of the number of the passes \cite{The_art_of_computer_programming_volume_3:_(2nd_ed.)_sorting_and_searching}. However, due to the intractability of asymptotic approximation, experimental study rather than theoretical analysis has always been adopted.
    
    This paper relates bubble sort to birthday problem. In section \ref{section_birthday} and section \ref{section_bubble}, birthday problem and bubble sort are defined formally. In section \ref{section_Rayleigh}, based on a generalization of Poisson limit theorem for dissociated random variables stated as theorem \ref{theorem_generalization_of_poisson}, an intuitive explanation indicates how Rayleigh distribution naturally arises in birthday problem and bubble sort in theorem \ref{theorem_collision_convergence} and theorem \ref{theorem_pass_convergence}, which proves that $\displaystyle \frac{n - P_{n}}{\sqrt{n}}$ converges in distribution to the standard Rayleigh distribution. In section \ref{section_estimation}, bounding bubble sort by birthday problem leads to a novel direction for analysing bubble sort. In section \ref{section_distribution} and section \ref{section_characteristics}, asymptotic distributions and statistical characteristics of bubble sort and birthday problem are presented. In section \ref{section_optimizations}, it is discovered that some common optimizations of bubble sort could lead to average performance degradation, which is the original motivation of this paper.

\section{Birthday problem}\label{section_birthday}

    Birthday problem asks for the probability that $m + 1$ persons have distinct birthdays in a year of $n$ days, which is determined by the distribution of the number of persons required to generate a birthday collision \cite{Analytic_Combinatorics}. Its isomorphic form in cryptography is known as birthday attack, which concentrates on the number of attacks required to generate a hash collision \cite{Introduction_to_Modern_Cryptography}. Moreover, birthday problem for $m = 22, n = 365$ is often mentioned in textbooks such as \cite{Probability:_Theory_and_Examples, The_art_of_computer_programming_volume_3:_(2nd_ed.)_sorting_and_searching}, because the probability that $23$ persons have distinct birthdays in a year of $365$ days is close to $\displaystyle \frac{1}{2}$.
    
    Note that this paper restricts birthday problem to the case that birthdays of persons are independent and identically uniformly distributed, which is consistent with definition \ref{definition_birthday}.
    \begin{definition}\label{definition_birthday}
        $U_{n, i}$ represents the birthday of the $i$-th person, where $n$ is the number of days in a year. Formally, $U_{n, 1}, U_{n, 2}, U_{n, 3}, \cdots$ are independent and identically uniformly distributed on $\left\{ 0, 1, \cdots, n - 1 \right\}$. Note that this definition is similar to definition \ref{definition_inversion}.
    \end{definition}
    \begin{definition}\label{definition_collision}
        $C_{n}$ represents the number of persons at the first collision, where $n$ is the number of days in a year. Formally, $C_{n} = \min \left\{ i + 1 \in \mathbf{N}^{+} \mid U_{n, i + 1} \in \left\{ U_{n, 1}, U_{n, 2}, \cdots, U_{n, i} \right\} \right\}$.
    \end{definition}
    \begin{lemma}\label{lemma_birthday_problem}
        $C_{n} > m + 1$ if and only if $\left| \left\{ U_{n, 1}, U_{n, 2}, \cdots, U_{n, m + 1} \right\} \right| = m + 1$.
    \end{lemma}
    \begin{proof}
        Use definition \ref{definition_collision}.
    \end{proof}
    \begin{theorem}\label{theorem_birthday_distinction}
        The persons at the first collision can be related to the distinction of the number of the birthdays as
        \begin{equation*}
            \mathbb{P} \left\{ C_{n} > m + 1 \right\} = \mathbb{P} \left\{ \left| U_{n, 1}, U_{n, 2}, \cdots, U_{n, m + 1} \right| = m + 1 \right\}
        \end{equation*}
    \end{theorem}
    \begin{proof}
        Use lemma \ref{lemma_birthday_problem}. Note that this theorem is isomorphic to theorem \ref{theorem_inversion_distinction}.
    \end{proof}
    \begin{theorem}\label{theorem_collision}
        The distribution of the number of persons at the first collision is
        \begin{equation*}
            \mathbb{P} \left\{ C_{n} > m + 1 \right\} = \left( 1 - \frac{1}{n} \right) \left( 1 - \frac{2}{n} \right) \cdots \left( 1 - \frac{m}{n} \right)
        \end{equation*}
    \end{theorem}
    \begin{proof}
        Use definition \ref{definition_collision} and lemma \ref{lemma_birthday_problem}. Note that this theorem is similar to theorem \ref{theorem_pass}.
    \end{proof}
    
\section{Bubble sort}\label{section_bubble}

    Bubble sort is composed of passes which traverse the sequence to swap adjacent elements out of order \cite{Introduction_to_Algorithms_Third_Edition}, as algorithm \ref{algorithm_bubble_sort} displays. Its static single-assignment form is presented by algorithm \ref{algorithm_static_single_assignment_form_of_bubble_sort} with definition \ref{definition_swap}. Moreover, terminating if no swaps are in a pass is always regarded as an optimization of bubble sort, which is described in algorithm \ref{algorithm_optimization_of_bubble_sort}.

    Note that this paper restricts bubble sort to the cases where the input is a uniformly random permutation of distinct numbers, which is consistent with definition \ref{definition_inversion}.

    \begin{definition}\label{definition_swap}
        $\varsigma_{j}$ represents the swap of the $j$-th element and the $\left( j + 1 \right)$-th element of a sequence. Formally, $\varsigma_{j} \left( a_{1}, a_{2}, \cdots, a_{n} \right) = \left( a'_{1}, a'_{2}, \cdots, a'_{n} \right)$ where $a'_{j} = a_{j + 1}, a'_{j + 1} = a_{j}$ and $a'_{k} = a_{k}$ for $k \not \in \left\{ j, j + 1 \right\}$.
    \end{definition}
    \begin{definition}
        $\left( a_{1}^{\left( i, n - i \right)}, a_{2}^{\left( i, n - i \right)}, \cdots, a_{n}^{\left( i, n - i \right)} \right)$ represents the value of $\left( a_{1}, a_{2}, \cdots, a_{n} \right)$ at the end of the $i$-th pass. Formally, $\left( a_{1}^{\left( i, j \right)}, a_{2}^{\left( i, j \right)}, \cdots, a_{n}^{\left( i, j \right)} \right)$ is recursively defined in algorithm \ref{algorithm_static_single_assignment_form_of_bubble_sort}.
    \end{definition}
    \begin{definition}\label{definition_pass}
        $P_{n}$ represents the number of passes required to bubble sort $\left( a_{1}, a_{2}, \cdots, a_{n} \right)$, as algorithm \ref{algorithm_optimization_of_bubble_sort} terminates at the end of the $P_{n}$-th pass. Formally, $P_{n} = \min \left\{ i + 1 \in  \mathbf{N}^{+} \mid a_{1}^{\left( i, n - i \right)} \leqslant a_{2}^{\left( i, n - i \right)} \leqslant \cdots \leqslant a_{n}^{\left( i, n - i \right)} \right\}$.
    \end{definition}
    \begin{definition}\label{definition_inversion}
        $V_{n, i}$ represents the number of inversions whose smaller component is the $i$-th smallest element. Formally, $V_{n, 1}, V_{n, 2}, \cdots, V_{n, n}$ are independent and $V_{n, i}$ is uniformly distributed on $\left\{ 0, 1, \cdots, n - i \right\}$. Note that this definition is similar to definition \ref{definition_birthday}.
    \end{definition}
    \begin{lemma}\label{lemma_bubble_sort}
       Passes of bubble sort are determined by the inversions as $P_{n} = \max \left\{ V_{n, 1} + 1, V_{n, 2} + 1, \cdots, V_{n, n} + 1 \right\}$.
    \end{lemma}
    \begin{proof}
        See \cite{The_art_of_computer_programming_volume_3:_(2nd_ed.)_sorting_and_searching}.
    \end{proof}
    \begin{theorem}\label{theorem_inversion_distinction}
        The passes required can be related to the distinction of the number of the inversions as
        \begin{equation*}
            \mathbb{P} \left\{ P_{n} \leqslant n - m \right\} = \mathbb{P} \left\{ \left| V_{n, 1}, V_{n, 2}, \cdots, V_{n, m + 1} \right| = m + 1 \right\}
        \end{equation*}
    \end{theorem}
    \begin{proof}
        Use definition \ref{definition_inversion} and lemma \ref{lemma_bubble_sort}. Note that this theorem is isomorphic to theorem \ref{theorem_birthday_distinction}.
    \end{proof}
    \begin{theorem}\label{theorem_pass}
        The distribution of the number of passes required is
        \begin{equation*}
            \mathbb{P} \left\{ P_{n} \leqslant n - m \right\} = \left( 1 - \frac{1}{n - m + 1} \right) \left( 1 - \frac{2}{n - m + 2} \right) \cdots \left( 1 - \frac{m}{n - m + m} \right)
        \end{equation*}
    \end{theorem}
    \begin{proof}
        Use definition \ref{definition_inversion} and lemma \ref{lemma_bubble_sort}. Note that this theorem is similar to \ref{theorem_collision}.
    \end{proof}
    
\section{Rayleigh distribution in birthday problem and bubble sort}\label{section_Rayleigh}
    Poisson distribution, exponential distribution and Rayleigh distribution are fundamental and elementary in probability theory. For a Poisson process, as lemma \ref{lemma_Poisson_process} shows, Poisson distribution represents the number of the events, and exponential distribution represents the distance between consecutive events. Moreover, Rayleigh distribution is related to exponential distribution by lemma \ref{lemma_relate_Rayleigh_to_exponential}.

    Note that theorem \ref{theorem_generalization_of_poisson} is a generalization of Poisson limit theorem for dissociated random variables, then theorem \ref{theorem_collision_convergence} and theorem \ref{theorem_pass_convergence} offer an intuitive explanation indicating how Rayleigh distribution naturally arises in birthday problem and bubble sort.

    \begin{definition}\label{definition_Kronecker_delta}
        $\delta$ represents Kronecker delta. Formally, $\delta \left( X, Y \right) = 1$ for $X = Y$, whereas $\delta \left( X, Y \right) = 0$ for $X \neq Y$.
    \end{definition}
    \begin{definition}\label{definition_Poisson}
        $\mathrm{Poi} \left( \lambda \right)$ represents Poisson distribution, and $\mathrm{Poi} \left( 1 \right)$ represents the standard Poisson distribution. Formally, $W$ follows $\mathrm{Poi} \left( \lambda \right)$ if and only if $\displaystyle \mathbb{P} \left\{ W = w \right\} = \frac{\lambda^{w}}{w !} \mathrm{e}^{- \lambda}$ for $w \in \mathbf{N}$.
    \end{definition}
    \begin{definition}\label{definition_exponential}
        $\mathrm{Exp} \left( \lambda \right)$ represents exponential distribution, and $\mathrm{Exp} \left( 1 \right)$ represents the standard exponential distribution. Formally, $W$ follows $\mathrm{Exp} \left( \lambda \right)$ if and only if $\displaystyle \mathbb{P} \left\{ W > w \right\} = \mathrm{e}^{- \lambda w}$ for $w \in \mathbf{R}^{+}$.
    \end{definition}
    \begin{definition}\label{definition_Rayleigh}
        $\mathrm{Ray} \left( \lambda \right)$ represents Rayleigh distribution, and $\mathrm{Ray} \left( 1 \right)$ represents the standard Rayleigh distribution. Formally, $W$ follows $\mathrm{Ray} \left( \lambda \right)$ if and only if $\displaystyle \mathbb{P} \left\{ W > w \right\} = \mathrm{e}^{- \frac{w^{2}}{2 \lambda^{2}}}$ for $w \in \mathbf{R}^{+}$.
    \end{definition}
    \begin{definition}\label{definition_Poisson_process}
        Poisson process represents a counting process with independent increments and stationary increments that the number of events follows Poisson distribution. Formally, $\left\{ Q \left( t \right) \mid t \geqslant 0 \right\}$ is a Poisson process with the parameter $\lambda$ if and only if $\displaystyle \mathbb{P} \left\{ Q \left( s + t \right) - Q \left( s \right) = k \right\} = \frac{ \left( \lambda t \right)^{k}}{k !} \mathrm{e}^{- \lambda t}$ for $s \geqslant 0, t \geqslant 0, k \in \mathbf{N}$ and the elements of $\left\{ Q \left( t_{i + 1} \right) - Q \left( t_{i} \right) \mid i \in \mathbf{N}_{+} \right\}$ are independent for $0 \leqslant t_{1} \leqslant t_{2} \leqslant t_{3} \leqslant \cdots$.
    \end{definition}
    \begin{definition}\label{definition_dissociation}
        Dissociation represents a generalization of independence for indicator random variables. Formally, the elements of $\left\{ D_{\left\{ i, j \right\}} \mid \left\{ i, j \right\} \in \mathscr{S} \right\}$ are dissociated if and only if $\left\{ D_{\left\{ i, j \right\}} \mid \left\{ i, j \right\} \in \mathscr{A} \right\}$ and $\left\{ D_{\left\{ i, j \right\}} \mid \left\{ i, j \right\} \in \mathscr{B} \right\}$ are independent for $\mathscr{A} \cup \mathscr{B} \subset \mathscr{S}$ and $\displaystyle \left( \bigcup_{\left\{ i, j \right\} \in \mathscr{A}} \left\{ i , j \right\} \right) \cap \left( \bigcup_{\left\{ i, j \right\} \in \mathscr{B}} \left\{ i , j \right\} \right) = \varnothing$.
    \end{definition}
    \begin{lemma}\label{lemma_Poisson_convergencence}
        If $W_{n}$ converges in distribution to $\mathrm{Poi} \left( \lambda \right)$, then $\mathbb{E} \left( W_{n} \right) \to \lambda$ and $\mathbb{V} \left( W_{n} \right) \to \lambda$.
    \end{lemma}
    \begin{proof}
        Use definition \ref{definition_Poisson}.
    \end{proof}
    \begin{lemma}\label{lemma_Rayleigh_convergencence}
        If $W_{n}$ converges in distribution to $\mathrm{Ray} \left( 1 \right)$, then $\displaystyle \mathbb{E} \left( W_{n}^{k} \right) \to \sqrt{2}^{k} \ \Gamma \left( \frac{k}{2} + 1 \right)$ for $k \in \mathbf{N}$ and $\varphi_{W_{n}} \left( t \right) \to \sqrt{2 \pi} \ \mathrm{e}^{- \frac{t^{2}}{2}} \ \mathrm{i} t \ \Phi \left( \mathrm{i} t \right) + 1$. Note that this lemma is mentioned in section \ref{section_characteristics}.
    \end{lemma}
    \begin{proof}
        Use definition \ref{definition_Rayleigh}.
    \end{proof}
    \begin{lemma}\label{lemma_relate_Rayleigh_to_exponential}
        $\displaystyle \frac{W_{n}}{\sqrt{n}}$ converges in distribution to $\mathrm{Ray} \left( 1 \right)$ if and only if $\displaystyle \frac{W_{n} \left( W_{n} + 1 \right)}{2 n}$ converges in distribution to $\mathrm{Exp} \left( 1 \right)$.
    \end{lemma}
    \begin{proof}
        Use definition definition \ref{definition_exponential} and definition \ref{definition_Rayleigh}. Note that $W$ follows $\mathrm{Ray} \left( 1 \right)$ if and only if $\displaystyle \frac{W^{2}}{2}$ follows $\mathrm{Exp} \left( 1 \right)$.
    \end{proof}
    \begin{lemma}\label{lemma_Poisson_process}
        If $\left\{ Q \left( t \right) \mid t \geqslant 0 \right\}$ is a Poisson process with the parameter $\lambda$, then the number of the events occurring between $s$ and $s + t$ is $Q \left( s + t \right) - Q \left( s \right)$ which follows $\mathrm{Poi} \left( \lambda t \right)$, and the distance between consecutive events is $\min \left\{ t \geqslant 0 \mid Q \left( s + t \right) > Q \left( s \right) \right\}$ which follows $\mathrm{Exp} \left( \lambda \right)$.
    \end{lemma}
    \begin{proof}
        Use definition \ref{definition_Poisson}, definition \ref{definition_exponential} and definition \ref{definition_Poisson_process}.
    \end{proof}
    \begin{lemma}\label{lemma_birthday_dissociation}
        The elements of $\left\{ \delta \left( U_{n, i}, U_{n, j} \right) \mid i < j \right\}$ are dissociated and $\left\{ 0, 1 \right\}$-valued.
    \end{lemma}
    \begin{proof}
        Use definition \ref{definition_birthday} and definition \ref{definition_dissociation}. Note that this lemma is isomorphic to lemma \ref{lemma_inversion_dissociation}.
    \end{proof}
    \begin{lemma}\label{lemma_inversion_dissociation}
        The elements of $\left\{ \delta \left( V_{n, i}, V_{n, j} \right) \mid i < j \right\}$ are dissociated and $\left\{ 0, 1 \right\}$-valued.
    \end{lemma}
    \begin{proof}
        Use definition \ref{definition_inversion} and definition \ref{definition_dissociation}. Note that this lemma is isomorphic to lemma \ref{lemma_birthday_dissociation}.
    \end{proof}
    \begin{theorem}\label{theorem_generalization_of_poisson}
        For $\displaystyle \mathbb{E} \left( \sum_{\left\{ i, j \right\} \in \mathscr{S}_{n}} D_{n, \left\{ i, j \right\}} \right) \to \lambda$ and $\displaystyle \mathbb{V} \left( \sum_{\left\{ i, j \right\} \in \mathscr{S}_{n}} D_{n, \left\{ i, j \right\}} \right) \to \lambda$, if the elements of $\left\{ D_{n, \left\{ i, j \right\}} \mid \left\{ i, j \right\} \in \mathscr{S}_{n} \right\}$ are dissociated and $\left\{ 0, 1 \right\}$-valued with $\displaystyle \left| T_{n} \right| \sum_{\left\{ i, j \right\} \in \mathscr{S}_{n}} \left( \mathbb{E} \left( D_{n, \left\{ i, j \right\}} \right) \right)^{2} \to 0$ and $\displaystyle \sum_{i \in T_{n}} \sum_{j \in T_{n} \setminus \left\{ i \right\}} \sum_{k \in T_{n} \setminus \left\{ i, j \right\}} \mathbb{E} \left( D_{n, \left\{ i, j \right\}} D_{n, \left\{ i, k \right\}} \right) \to 0$ where $\displaystyle T_{n} = \bigcup_{\left\{ i, j \right\} \in \mathscr{S}_{n}} \left\{ i, j \right\}$, then $\displaystyle \sum_{\left\{ i, j \right\} \in \mathscr{S}_{n}} D_{n, \left\{ i, j \right\}}$ converges in distribution to $\mathrm{Poi} \left( \lambda \right)$.
    \end{theorem}
    \begin{proof}
        See appendix \ref{proof_generalization_of_poisson_theorem}. Note that this theorem is consistent with lemma \ref{lemma_Poisson_convergencence}.
    \end{proof}
    \begin{theorem}\label{theorem_collision_count}
        $\displaystyle \sum_{1 \leqslant i < j \leqslant m + 1} \delta \left( U_{n, i}, U_{n, j} \right)$ converges in distribution to $\mathrm{Poi} \left( \lambda \right)$ where $\displaystyle \frac{m \left( m + 1 \right)}{2 n} \to \lambda$.
    \end{theorem}
    \begin{proof}
        Use lemma \ref{lemma_birthday_dissociation} and theorem \ref{theorem_generalization_of_poisson}. Note that this theorem is isomorphic to theorem \ref{theorem_pass_count}.
    \end{proof}
    \begin{theorem}\label{theorem_pass_count}
        $\displaystyle \sum_{1 \leqslant i < j \leqslant m + 1} \delta \left( V_{n, i}, V_{n, j} \right)$ converges in distribution to $\mathrm{Poi} \left( \lambda \right)$ where $\displaystyle \frac{m \left( m + 1 \right)}{2 n} \to \lambda$.
    \end{theorem}
    \begin{proof}
        Use lemma \ref{lemma_inversion_dissociation} and theorem \ref{theorem_generalization_of_poisson}. Note that this theorem is isomorphic to theorem \ref{theorem_collision_count}.
    \end{proof}
    \begin{theorem}\label{theorem_collision_convergence}
        $\displaystyle \frac{C_{n} - 1}{\sqrt{n}}$ converges in distribution to the standard Rayleigh distribution.
    \end{theorem}
    \begin{proof}
        Inspired by lemma \ref{lemma_Poisson_process}, using theorem \ref{theorem_birthday_distinction} and theorem \ref{theorem_collision_count} gives that $\displaystyle \frac{\left( C_{n} - 1 \right) \left( C_{n} - 1 + 1 \right)}{2 n}$ converges in distribution to the standard exponential distribution, then lemma \ref{lemma_relate_Rayleigh_to_exponential} naturally leads to the standard Rayleigh distribution. Note that this theorem is isomorphic to theorem \ref{theorem_pass_convergence}.
    \end{proof}
    \begin{theorem}\label{theorem_pass_convergence}
        $\displaystyle \frac{n - P_{n}}{\sqrt{n}}$ converges in distribution to the standard Rayleigh distribution.
    \end{theorem}
    \begin{proof}
        Inspired by lemma \ref{lemma_Poisson_process}, using theorem \ref{theorem_inversion_distinction} and theorem \ref{theorem_pass_count} gives that $\displaystyle \frac{\left( n - P_{n} \right) \left( n - P_{n} + 1 \right)}{2 n}$ converges in distribution to the standard exponential distribution, then lemma \ref{lemma_relate_Rayleigh_to_exponential} naturally leads to the standard Rayleigh distribution. Note that this theorem is isomorphic to theorem \ref{theorem_collision_convergence}.
    \end{proof}
    
\section{Estimation of bubble sort with birthday problem}\label{section_estimation}

    As a well-known classic problem in probability theory widely applied in cryptography, birthday problem has been extensively studied. However, bubble sort does not receive much attention from researchers due to its inefficiency, although it is a popular sorting algorithm leading to some interesting theoretical problems in computer science.

    Note that lemma \ref{lemma_bound} bounds bubble sort by birthday problem, then estimating bubble sort with birthday problem presents a novel direction for analysing bubble sort.

    \begin{lemma}\label{lemma_bound}
        Bubble sort can be bounded by birthday problem as
        \begin{equation*}
            \mathbb{P} \left\{ C_{n - \left( m - 1 \right)} > m + 1 \right\} \leqslant \mathbb{P} \left\{ P_{n} \leqslant n - m \right\} \leqslant \mathbb{P} \left\{ C_{n} > m + 1 \right\}
        \end{equation*}
    \end{lemma}
    \begin{proof}
        Use theorem \ref{theorem_collision} and theorem \ref{theorem_pass}. Note that $\displaystyle \frac{C_{n} - 1}{\sqrt{n}}$ follows the standard Rayleigh distribution \cite{Probability:_Theory_and_Examples}, then using this theorem can simply give that $\displaystyle \frac{n - P_{n}}{\sqrt{n}}$ follows the standard Rayleigh distribution, which is stated as theorem \ref{theorem_pass_convergence}.
    \end{proof}
    \begin{lemma}\label{lemma_birthday_series}
        The distribution of the number of persons at the first collision is
        \begin{equation*}
            \mathbb{P} \left\{ C_{n} > m + 1 \right\} = \exp \sum_{k \in \mathrm{N}_{+}} \frac{1^{k} + 2^{k} + \cdots + m^{k}}{- k n^{k}}
        \end{equation*}
    \end{lemma}
    \begin{proof}
        Use theorem \ref{theorem_collision} and Maclaurin series of $\ln \left( 1 - x \right)$.
    \end{proof}
    \begin{lemma}\label{lemma_bubble_series}
        The distribution of the number of passes required is
        \begin{equation*}
            \mathbb{P} \left\{ P_{n} \leqslant n - m \right\} = \exp \sum_{k \in \mathrm{N}_{+}} \frac{1^{k} + 2^{k} + \cdots + m^{k}}{\left( - 1 \right)^{k} k \left( n - m \right)^{k}}
        \end{equation*}
    \end{lemma}
    \begin{proof}
        Use theorem \ref{theorem_pass} and Maclaurin series of $\ln \left( 1 + x \right)$.
    \end{proof}
    \begin{theorem}
        The distribution of the number of persons at the first collision is
        \begin{equation*}
            \mathbb{P} \left\{ C_{n} > m + 1 \right\} = \exp \frac{m \left( m + 1 \right)}{\displaystyle - 2 \left( n - \frac{m}{3} \right)} \pm \Theta \left( \frac{m^{4}}{n^{3}} \right)
        \end{equation*}
    \end{theorem}
    \begin{proof}
        Use lemma \ref{lemma_birthday_series}.
    \end{proof}
    \begin{theorem}
        The distribution of the number of passes required is
        \begin{equation*}
            \mathbb{P} \left\{ P_{n} \leqslant n - m \right\} = \exp \frac{m \left( m + 1 \right)}{\displaystyle - 2 \left( n - \frac{2 m}{3} \right)} \pm \Theta \left( \frac{m^{4}}{n^{3}} \right)
        \end{equation*}
    \end{theorem}
    \begin{proof}
        Use lemma \ref{lemma_bubble_series}.
    \end{proof}
    \begin{theorem}\label{theorem_common_estimation}
        The relative error of estimating $\mathbb{P} \left\{ P_{n} \leqslant n - m \right\}$ with $\mathbb{P} \left\{ C_{n} > m + 1 \right\}$ is
        \begin{equation*}
            \frac{\mathbb{P} \left\{ C_{n} > m + 1 \right\}}{\mathbb{P} \left\{ P_{n} \leqslant n - m \right\}} - 100 \% = \frac{\left( m - 1 \right) m \left( m + 1 \right)}{\displaystyle 6 \left( n - \frac{m}{2} \right)^{2}} \pm \Theta \left( \frac{m^{5}}{n^{4}} \right)
        \end{equation*}
    \end{theorem}
    \begin{proof}
        Use lemma \ref{lemma_birthday_series} and lemma \ref{lemma_bubble_series}.
    \end{proof}
    \begin{theorem}\label{theorem_optimal_estimation}
        The relative error of estimating $\mathbb{P} \left\{ P_{n} \leqslant n - m \right\}$ with $\mathbb{P} \left\{ C_{n - \frac{m - 1}{3}} > m + 1 \right\}$ is
        \begin{equation*}
            \frac{\mathbb{P} \left\{ C_{n - \frac{m - 1}{3}} > m + 1 \right\}}{\mathbb{P} \left\{ P_{n} \leqslant n - m \right\}} - 100 \% = \frac{\left( m - 1 \right) m \left( m + 1 \right) \left( m + 2 \right) \left( 2 m + 1 \right)}{\displaystyle - 270 \left( n - \frac{4 m - 1}{6} \right)^{4}} \pm \Theta \left( \frac{m^{7}}{n^{6}} \right)
        \end{equation*}
    \end{theorem}
    \begin{proof}
        Use lemma \ref{lemma_birthday_series} and lemma \ref{lemma_bubble_series}.
    \end{proof}
    \begin{theorem}\label{theorem_estimation}
        The optimal argument $k$ for estimating $\mathbb{P} \left\{ P_{n} \leqslant n - m \right\}$ with $\mathbb{P} \left\{ C_{n - k} > m + 1 \right\}$ is
        \begin{equation*}
            \underset{0 \leqslant k < m}{\mathrm{argmin}} \left| \frac{\mathbb{P} \left\{ C_{n - k} > m + 1 \right\}}{\mathbb{P} \left\{ P_{n} \leqslant n - m \right\}} - 100 \% \right| = \frac{m - 1}{3} \pm \Theta \left( \frac{m^{3}}{n^{2}} \right)
        \end{equation*}
    \end{theorem}
    \begin{proof}
        Use lemma \ref{lemma_birthday_series} and lemma \ref{lemma_bubble_series}.
    \end{proof}
    \begin{example}
        See theorem \ref{theorem_common_estimation} and theorem \ref{theorem_optimal_estimation}, then substituting $m = 22, n = 365$ gives
        \begin{equation*}
            \begin{aligned}
                \mathbb{P} \left\{ C_{358} >          23 \right\} & \approx 0.4857834 \\
                \mathbb{P} \left\{ P_{365} \leqslant 343 \right\} & \approx 0.4857848 \\
                \mathbb{P} \left\{ C_{365} >          23 \right\} & \approx 0.4927028
            \end{aligned}
        \end{equation*}
    \end{example}

\section{Asymptotic distributions of bubble sort and birthday problem}\label{section_distribution}

    The factorial function is pivotal for analysing the distributions of bubble sort and birthday problem as theorem \ref{theorem_collision} and theorem \ref{theorem_pass} present, which leads to lemma \ref{lemma_bubble_fractorial} and lemma \ref{lemma_birthday_fractorial}. Then the cumulative distribution functions and the probability mass functions are presented.
    
    \begin{definition}\label{definition_bubble_standard}
        $X_{n}$ represents the value of an affine function at $P_{n}$, inspired by theorem \ref{theorem_pass_convergence}. Formally, $\displaystyle X_{n} = \frac{n - P_{n}}{\sqrt{n}}$.
    \end{definition} \begin{definition}\label{definition_birthday_standard}
        $Z_{n}$ represents the value of an affine function at $C_{n}$, inspired by theorem \ref{theorem_collision_convergence}. Formally, $\displaystyle Z_{n} = \frac{C_{n} - 1}{\sqrt{n}}$.
    \end{definition}
    \begin{definition}\label{definition_varrho}
        $\varrho_{n}$ represents an extension of the distribution of $X_{n}$, inspired by lemma \ref{lemma_bubble_fractorial}. Formally, $\displaystyle \varrho_{n} \left( x \right) = \frac{\Gamma \left( n - x \sqrt{n} + 1 \right)}{\Gamma \left( n + 1 \right)} \left( n - x \sqrt{n} \right)^{x \sqrt{n}}$.
    \end{definition}
    \begin{lemma}\label{lemma_Stirling}
        For $x \geqslant 1$, there holds
        \begin{equation*}
            \Gamma \left( x + 1 \right) = \sqrt{2 \pi x} \left( \frac{x}{\mathrm{e}} \right)^{x} \exp \left( \frac{1}{12 x} \pm \Theta \left( \frac{1}{x^{3}} \right) \right)
        \end{equation*}
    \end{lemma}
    \begin{proof}
        Use Stirling series.
    \end{proof}
    \begin{lemma}\label{lemma_bubble_fractorial}
        For $n - x \sqrt{n} \in \left\{ 1, 2, \cdots, n \right\}$, the distribution of $X_{n}$ is
        \begin{equation*}
            \mathbb{P} \left\{ X_{n} \geqslant x \right\} = \frac{\left( n - x \sqrt{n} \right) !}{n !}  \left( n - x \sqrt{n} \right)^{x \sqrt{n}}
        \end{equation*}
    \end{lemma}
    \begin{proof}
        Use theorem \ref{theorem_pass} and definition \ref{definition_bubble_standard}.
    \end{proof}
    \begin{lemma}\label{lemma_birthday_fractorial}
        For $z \sqrt{n} \in \left\{ 1, 2, \cdots, n \right\}$, the distribution of $Z_{n}$ is
        \begin{equation*}
            \mathbb{P} \left\{ Z_{n} \geqslant z \right\} = \frac{n !}{\left( n - z \sqrt{n} \right) !} \left( \frac{1}{n} \right)^{z \sqrt{n}} 
        \end{equation*}
    \end{lemma}
    \begin{proof}
        Use theorem \ref{theorem_collision} and definition \ref{definition_birthday_standard}.
    \end{proof}
    \begin{theorem}\label{theorem_varrho}
        For $x \geqslant 0$, there holds
        \begin{equation*}
            \varrho_{n} \left( x \right) = \mathrm{e}^{-\frac{x^{2}}{2}} \exp \left( - \frac{2 x^{3} + 3 x}{6 \sqrt{n}} - \frac{x^{4} + x^{2}}{4 n} - \frac{12 x^{5} + 10 x^{3} - 5 x}{60 n \sqrt{n}} - \frac{4 x^{6} + 3 x^{4} - 2 x^{2}}{24 n^{2}} \pm \Theta \left( \frac{x^{7}}{n^{2} \sqrt{n}} \right) \right)
        \end{equation*}
    \end{theorem}
    \begin{proof}
        Use definition \ref{definition_varrho} and lemma \ref{lemma_Stirling}.
    \end{proof}
    \begin{theorem}\label{theorem_bubble_CDF}
        For $n - x \sqrt{n} \in \left\{ 1, 2, \cdots, n \right\}$, the cumulative distribution function (CDF) of $X_{n}$ is
        \begin{equation*}
            F_{X_{n}} \left( x \right) = 1 - \mathrm{e}^{-\frac{x^{2}}{2}} \exp \left( - \frac{2 x^{3} + 9 x}{6 \sqrt{n}} \pm \Theta \left( \frac{x^{4}}{n} \right) \right)
        \end{equation*}
    \end{theorem}
    \begin{proof}
        Use theorem \ref{theorem_varrho} with $\displaystyle F_{X_{n}} \left( x \right) = 1 - \varrho_{n} \left( x + \frac{1}{\sqrt{n}} \right)$. Note that this theorem is consistent with theorem \ref{theorem_pass_convergence}.
    \end{proof}
    \begin{theorem}\label{theorem_bubble_PMF}
        For $n - x \sqrt{n} \in \left\{ 1, 2, \cdots, n \right\}$, the probability mass function (PMF) of $X_{n}$ is
        \begin{equation*}
            f_{X_{n}} \left( x \right) = \frac{\mathrm{e}^{- \frac{x^{2}}{2}} x}{\sqrt{n}} \exp \left( - \frac{x^{4} - 3}{3 x \sqrt{n}} \pm \Theta \left( \frac{x^{4}}{n} \right) \right)
        \end{equation*}
    \end{theorem}
    \begin{proof}
        Use theorem \ref{theorem_varrho} with $\displaystyle f_{X_{n}} \left( x \right) = \varrho_{n} \left( x \right) - \varrho_{n} \left( x + \frac{1}{\sqrt{n}} \right)$.
    \end{proof}
    \begin{theorem}\label{theorem_birthday_CDF}
        For $z \sqrt{n} \in \left\{ 1, 2, \cdots, n \right\}$, the cumulative distribution function (CDF) of $Z_{n}$ is
        \begin{equation*}
            F_{Z_{n}} \left( z \right) = 1 - \mathrm{e}^{- \frac{z^{2}}{2}} \exp \left( - \frac{z^{3} + 3 z}{6 \sqrt{n}} \pm \Theta \left( \frac{z^{4}}{n} \right) \right)
        \end{equation*}
    \end{theorem}
    \begin{proof}
        Use lemma \ref{lemma_Stirling} and lemma \ref{lemma_birthday_fractorial}. Note that this theorem is consistent with theorem \ref{theorem_collision_convergence}.
    \end{proof}
    \begin{theorem}\label{theorem_birthday_PMF}
        For $z \sqrt{n} \in \left\{ 1, 2, \cdots, n \right\}$, the probability mass function (PMF) of $Z_{n}$ is
        \begin{equation*}
            f_{Z_{n}} \left( z \right) = \frac{\mathrm{e}^{- \frac{z^{2}}{2}} z}{\sqrt{n}} \exp \left( - \frac{z^{3} - 3 z}{6 \sqrt{n}} \pm \Theta \left( \frac{z^{4}}{n} \right) \right)
        \end{equation*}
    \end{theorem}
    \begin{proof}
        Use lemma \ref{lemma_Stirling} and lemma \ref{lemma_birthday_fractorial}.
    \end{proof}

\section{Statistical characteristics of bubble sort}\label{section_characteristics}

    Euler–Maclaurin summation formula is an effective approach for analysing the statistical characteristics of discrete random variables, such as expectation, variance, moment and characteristic function. Then these statistical characteristics of bubble sort are derived.

    \begin{definition}
        $\Phi$ represents the cumulative distribution function of the standard normal distribution. Formally, $\displaystyle \Phi \left( x \right) = \frac{1}{2} \left( \mathrm{erf} \left( \frac{x}{\sqrt{2}} \right) + 1 \right)$.
    \end{definition}
    \begin{lemma}\label{lemma_Euler_Maclaurin}
        For $\displaystyle 0 < \varepsilon < \frac{1}{6}$, the summation of $\varrho \left( x \right)$ is
        \begin{equation*}
            \sum_{n - x \sqrt{n} \in \left\{ 1, 2, \cdots, n \right\}} \varrho_{n} \left( x \right) = \sqrt{n} \int_{0}^{n^{\varepsilon}} \varrho_{n} \left( x \right) \mathrm{d} x + \frac{\varrho_{n} \left( 0 \right)}{2} - \frac{\varrho_{n}' \left( 0 \right)}{12 \sqrt{n}} + \frac{\varrho_{n}''' \left( 0 \right)}{720 n \sqrt{n}} \pm \Theta \left( \frac{1}{n^{2}} \right)
        \end{equation*}
    \end{lemma}
    \begin{proof}
        Use Euler–Maclaurin summation formula.
    \end{proof}
    \begin{lemma}\label{lemma_bubble_moment}
        For $k \in \mathbf{N}$, the $k$-th moment of $X_{n}$ is
        \begin{equation*}
            \mathbb{E} \left( X_{n}^{k} \right) = \left( - \frac{1}{\sqrt{n}} \right)^{k} + \sum_{n - x \sqrt{n} \in \left\{ 1, 2, \cdots, n \right\}} \left( x^{k} - \left( x - \frac{1}{\sqrt{n}} \right)^{k} \right) \varrho_{n} \left( x \right)
        \end{equation*}
    \end{lemma}
    \begin{proof}
        Use definition \ref{definition_varrho}, lemma \ref{lemma_bubble_fractorial} and Abel summation transformation.
    \end{proof}
    \begin{theorem}
        For $k \in \mathbf{N}$, the $k$-th moment of $X_{n}$ is
        \begin{equation*}
            \mathbb{E} \left( X_{n}^{k} \right) = \sqrt{2}^{k} \left( \Gamma \left( \frac{k}{2} + 1 \right) - \frac{\sqrt{2} k \left( k + 4 \right)}{6 \sqrt{n}} \Gamma \left( \frac{k + 1}{2} \right) \pm \Theta \left( \frac{1}{n} \right) \right)
        \end{equation*}
    \end{theorem}
    \begin{proof}
        Use lemma \ref{lemma_bubble_moment} and lemma \ref{lemma_Euler_Maclaurin}. Note that this theorem is consistent with lemma \ref{lemma_Rayleigh_convergencence}.
    \end{proof}
    \begin{theorem}
        The characteristic function of $X_{n}$ is
        \begin{equation*}
            \varphi_{X_{n}} \left( t \right) = \left( 1 - \frac{\left( 6 - t^{2} \right) \mathrm{i} t}{3 \sqrt{n}} \right) \sqrt{2 \pi} \ \mathrm{e}^{- \frac{t^{2}}{2}} \ \mathrm{i} t \ \Phi \left( \mathrm{i} t \right) + \left( 1 - \frac{\left( 5 - t^{2} \right) \mathrm{i} t}{3 \sqrt{n}} \right) \pm \Theta \left( \frac{1}{n} \right)
        \end{equation*}
    \end{theorem}
    \begin{proof}
        Use lemma \ref{lemma_Euler_Maclaurin}. Note that this theorem is consistent with lemma \ref{lemma_Rayleigh_convergencence}.
    \end{proof}
    \begin{theorem}\label{theorem_bubble_expectation}
        The expectation of $X_{n}$ is $\displaystyle \mathbb{E} \left( X_{n} \right) = \hat{\mathbb{E}}_{\frac{1}{n^{2} \sqrt{n}}} \left( X_{n} \right) \pm \Theta \left( \frac{1}{n^{2} \sqrt{n}} \right)$ where
        \begin{equation*}
            \hat{\mathbb{E}}_{\frac{1}{n^{2} \sqrt{n}}} \left( X_{n} \right) = \sqrt{\frac{\pi}{2}} - \frac{5}{3 \sqrt{n}} + \frac{11}{24 n} \sqrt{\frac{\pi}{2}} + \frac{4}{135 n \sqrt{n}} - \frac{71}{1152 n^{2}} \sqrt{\frac{\pi}{2}}
        \end{equation*}
    \end{theorem}
    \begin{proof}
        Use lemma \ref{lemma_Euler_Maclaurin}. Note that this theorem is consistent with Knuth's results presented in \cite{The_art_of_computer_programming_volume_1_(3rd_ed.):_fundamental_algorithm}.
    \end{proof}
    \begin{theorem}\label{theorem_bubble_square_expectation}
        The expectation of the square of $X_{n}$ is $\displaystyle \mathbb{E} \left( X_{n}^{2} \right) = \hat{\mathbb{E}}_{\frac{1}{n^{2} \sqrt{n}}} \left( X_{n}^{2} \right) \pm \Theta \left( \frac{1}{n^{2} \sqrt{n}} \right)$ where
        \begin{equation*}
            \hat{\mathbb{E}}_{\frac{1}{n^{2} \sqrt{n}}} \left( X_{n}^{2} \right) = 2 - 4 \sqrt{\frac{\pi}{2 n}} + \frac{5}{n} - \frac{5}{3 n} \sqrt{\frac{\pi}{2 n}} - \frac{4}{135 n^{2}}
        \end{equation*}
    \end{theorem}
    \begin{proof}
        Use lemma \ref{lemma_Euler_Maclaurin}.
    \end{proof}
    \begin{theorem}\label{theorem_bubble_variance}
        The variance of $X_{n}$ is $\displaystyle \mathbb{V} \left( X_{n} \right) = \hat{\mathbb{V}}_{\frac{1}{n^{2} \sqrt{n}}} \left( X_{n} \right) \pm \Theta \left( \frac{1}{n^{2} \sqrt{n}} \right)$ where
        \begin{equation*}
            \hat{\mathbb{V}}_{\frac{1}{n^{2} \sqrt{n}}} \left( X_{n} \right) = \frac{4 - \pi}{2} - \frac{2}{3} \sqrt{\frac{\pi}{2 n}} + \frac{160 - 33 \pi}{72 n} - \frac{107}{540 n} \sqrt{\frac{\pi}{2 n }} - \frac{1125 \pi - 1792}{25920 n^{2}}
        \end{equation*}
    \end{theorem}
    \begin{proof}
        Use theorem \ref{theorem_bubble_expectation} and theorem \ref{theorem_bubble_square_expectation}.
    \end{proof}
    \begin{example}
        See theorem \ref{theorem_bubble_expectation}, theorem \ref{theorem_bubble_square_expectation} and theorem \ref{theorem_bubble_variance}, then substituting $n = 10^{4}$ gives
        \begin{equation*}
            \begin{aligned}
                \mathbb{E}       \left( X_{10^{4}} \right)     & \approx 1.23670494307038 \\
                \hat{\mathbb{E}}_{\frac{1}{n^{2} \sqrt{n}}} \left( X_{10^{4}} \right)     & \approx 1.23670494307065 
            \end{aligned}
        \end{equation*}
        \begin{equation*}
            \begin{aligned}
                \mathbb{E}       \left( X_{10^{4}}^{2} \right) & \approx 1.950365345384   \\
                \hat{\mathbb{E}}_{\frac{1}{n^{2} \sqrt{n}}} \left( X_{10^{4}}^{2} \right) & \approx 1.950365345354
            \end{aligned}
        \end{equation*}
        \begin{equation*}
            \begin{aligned}
                \mathbb{V}       \left( X_{10^{4}} \right)     & \approx 0.4209262291695  \\
                \hat{\mathbb{V}}_{\frac{1}{n^{2} \sqrt{n}}} \left( X_{10^{4}} \right)     & \approx 0.4209262291679
            \end{aligned}
        \end{equation*}
    \end{example}

\section{Optimizations of bubble sort}\label{section_optimizations}

    Terminating if no swaps are in a pass is always regarded as an optimization of bubble sort, because it reduces the number of comparisons of adjacent elements. However, optimizations of bubble sort have to execute more other operations to determine if no swaps are in a pass, such as boolean assignments presented in algorithm \ref{algorithm_optimization_of_bubble_sort} and algorithm \ref{algorithm_variant_of_optimization_of_bubble_sort}. Therefore, simply regarding it as an optimization of bubble sort seems to be too arbitrary.

    Note that both one comparison of adjacent elements and one boolean assignment take $\Theta \left( 1 \right)$ time, in a model of computation which is similar to word random-access memory (RAM) model. Then the optimization of bubble sort (algorithm \ref{algorithm_optimization_of_bubble_sort}) must lead to average performance degradation due to theorem \ref{theorem_reduction} and theorem \ref{theorem_increasement}, whereas the variant of the optimization of bubble sort (algorithm \ref{algorithm_variant_of_optimization_of_bubble_sort}) might lead to average performance degradation due to theorem \ref{theorem_variant_reduction} and theorem \ref{theorem_variant_increasement}.

    \begin{lemma}\label{lemma_reduction}
        Compared to bubble sort (algorithm \ref{algorithm_bubble_sort}), the optimization of bubble sort (algorithm \ref{algorithm_optimization_of_bubble_sort}) reduces the number of comparisons of adjacent elements by $\displaystyle \frac{\left( n - P_{n} - 1 \right) \left( n - P_{n} \right)}{2}$.
    \end{lemma}
    \begin{proof}
        See algorithm \ref{algorithm_bubble_sort} and algorithm \ref{algorithm_optimization_of_bubble_sort}, then use definition \ref{definition_pass}.
    \end{proof}
    \begin{lemma}\label{lemma_variant_reduction}
        Compared to bubble sort (algorithm \ref{algorithm_bubble_sort}), the variant of the optimization of bubble sort (algorithm \ref{algorithm_variant_of_optimization_of_bubble_sort}) reduces the number of comparisons of adjacent elements by $\displaystyle \frac{\left( n - P_{n} - 1 \right) \left( n - P_{n} \right)}{2}$.
    \end{lemma}
    \begin{proof}
        See algorithm \ref{algorithm_bubble_sort} and algorithm \ref{algorithm_variant_of_optimization_of_bubble_sort}, then use definition \ref{definition_pass}.
    \end{proof}
    \begin{lemma}\label{lemma_increasement}
        Compared to bubble sort (algorithm \ref{algorithm_bubble_sort}), the optimization of bubble sort (algorithm \ref{algorithm_optimization_of_bubble_sort}) increases the number of boolean assignments by $\displaystyle P_{n} + \sum_{1 \leqslant i \leqslant n} V_{n, i}$.
    \end{lemma}
    \begin{proof}
        See algorithm \ref{algorithm_bubble_sort} and algorithm \ref{algorithm_optimization_of_bubble_sort}, then use definition \ref{definition_pass} and definition \ref{definition_inversion}.
    \end{proof}
    \begin{lemma}\label{lemma_variant_increasement}
        Compared to bubble sort (algorithm \ref{algorithm_bubble_sort}), the variant of the optimization of bubble sort (algorithm \ref{algorithm_optimization_of_bubble_sort}) increases the number of boolean assignments by $2 P_{n}$.
    \end{lemma}
    \begin{proof}
        See algorithm \ref{algorithm_bubble_sort} and algorithm \ref{algorithm_variant_of_optimization_of_bubble_sort}, then use definition \ref{definition_pass}.
    \end{proof}
    \begin{theorem}\label{theorem_reduction}
        Compared to bubble sort (algorithm \ref{algorithm_bubble_sort}), the optimization of bubble sort (algorithm \ref{algorithm_optimization_of_bubble_sort}) reduces the expectation of the number of comparisons of adjacent elements by
        \begin{equation*}
            n - \frac{5}{2} \sqrt{\frac{\pi n}{2}} + \frac{10}{3} - \frac{17}{16} \sqrt{\frac{\pi}{2 n}} - \frac{4}{135 n} \pm \Theta \left( \frac{1}{n \sqrt{n}} \right)
        \end{equation*}
    \end{theorem}
    \begin{proof}
        Use theorem \ref{theorem_bubble_expectation}, theorem \ref{theorem_bubble_square_expectation} and lemma \ref{lemma_reduction}.
    \end{proof}
    \begin{theorem}\label{theorem_variant_reduction}
        Compared to bubble sort (algorithm \ref{algorithm_bubble_sort}), the variant of the optimization of bubble sort (algorithm \ref{algorithm_variant_of_optimization_of_bubble_sort}) reduces the expectation of the number of comparisons of adjacent elements by
        \begin{equation*}
            n - \frac{5}{2} \sqrt{\frac{\pi n}{2}} + \frac{10}{3} - \frac{17}{16} \sqrt{\frac{\pi}{2 n}} - \frac{4}{135 n} \pm \Theta \left( \frac{1}{n \sqrt{n}} \right)
        \end{equation*}
    \end{theorem}
    \begin{proof}
        Use theorem \ref{theorem_bubble_expectation}, theorem \ref{theorem_bubble_square_expectation} and lemma \ref{lemma_variant_reduction}.
    \end{proof}
    \begin{theorem}\label{theorem_increasement}
        Compared to bubble sort (algorithm \ref{algorithm_bubble_sort}), the optimization of bubble sort (algorithm \ref{algorithm_optimization_of_bubble_sort}) increase the expectation of the number of boolean assignments by
        \begin{equation*}
            \frac{1}{2} n^{2} + \frac{1}{2} n - \sqrt{\frac{\pi n}{2}} + \frac{5}{3} - \frac{11}{24} \sqrt{\frac{\pi}{2 n}} - \frac{4}{135 n} \pm \Theta \left( \frac{1}{n \sqrt{n}} \right)
        \end{equation*}
    \end{theorem}
    \begin{proof}
        Use definition \ref{definition_inversion}, theorem \ref{theorem_bubble_expectation} and lemma \ref{lemma_increasement}.
    \end{proof}
    \begin{theorem}\label{theorem_variant_increasement}
        Compared to bubble sort (algorithm \ref{algorithm_bubble_sort}), the variant of the optimization of bubble sort (algorithm \ref{algorithm_variant_of_optimization_of_bubble_sort}) increase the expectation of the number of boolean assignments by
        \begin{equation*}
            2 n - 2 \sqrt{\frac{\pi n}{2}} + \frac{10}{3} - \frac{11}{12} \sqrt{\frac{\pi}{2 n}} - \frac{8}{135 n} \pm \Theta \left( \frac{1}{n \sqrt{n}} \right)
        \end{equation*}
    \end{theorem}
    \begin{proof}
        Use theorem \ref{theorem_bubble_expectation} and lemma \ref{lemma_variant_increasement}.
    \end{proof}

\section{Discussion}

    In analytic combinatorics, Kuba and Panholzer studied the distribution of inversions in labelled tree families and indicated how mixed Poisson-Rayleigh distributions naturally arise in several critical composition schemes \cite{On_moment_sequences_and_mixed_Poisson_distributions}. Inspired by their work, in section \ref{section_Rayleigh} and section \ref{section_estimation}, this paper provides a generalization of Poisson limit theorem for dissociated random variables to offer an intuitive explanation indicating how Rayleigh distribution naturally arises in birthday problem and bubble sort, and relating bubble sort to birthday problem presents a novel direction for analysing birthday problem and bubble sort. It is somewhat amazing that this relation has eluded discovery for so long.

    In computer science, Knuth studied the analysis of bubble sort and derived the expectation of its passes by the inversions \cite{The_art_of_computer_programming_volume_1_(3rd_ed.):_fundamental_algorithm, The_art_of_computer_programming_volume_3:_(2nd_ed.)_sorting_and_searching}. Based on his work, in section \ref{section_distribution} and section \ref{section_characteristics}, this paper develops this approach and derives several effective asymptotic approximations of bubble sort and birthday problem, including asymptotic distributions and statistical characteristics, which completely addresses the analysis of the average performance of bubble sort for distinct elements. And for analysing variations of bubble sort such as \cite{An_efficient_variation_of_bubble_sort}, some generalizations of this approach should be effective.

    Moreover, in computer engineering, nearly every description of bubble sort mentions the optimization that terminates if no swaps are made in a pass \cite{Bubble_sort:_an_archaeological_algorithmic_analysis}. However, in section \ref{section_optimizations}, this paper proves that some common implementations of this optimization such as \cite{Smart_Bubble_Sort:_A_Novel_and_Dynamic_Variant_of_Bubble_Sort_Algorithm, Experiment_Analysis_on_the_Bubble_Sort_Algorithm_and_Its_Improved_Algorithms} could lead to average performance degradation for sufficiently large $n$.

\newpage

\section{Acknowledgements}

    Sincere gratitude is expressed to Professor Ke Xu and Doctor Yicheng Pan for their encouragement and insightful suggestions throughout this research.

\bibliographystyle{cas-model2-names}
\bibliography{references}

\newpage

\appendix
\section{Algorithms}
    \subsection{Bubble sort}
    \begin{algorithm}[H]
        \caption{Bubble sort}
        \begin{algorithmic}[1]\label{algorithm_bubble_sort}
            \REQUIRE{$\left( a_{1}, a_{2}, \cdots, a_{n} \right)$}
            \ENSURE{a non-decreasing permutation of the input}
            \FOR{$i$ \textbf{from} $1$ \textbf{to} $n$}
                \FOR{$j$ \textbf{from} $1$ \textbf{to} $n - i$}
                    \IF{$a_{j} > a_{j + 1}$}
                        \STATE{swap $a_{j}, a_{j + 1}$}
                    \ENDIF
                \ENDFOR
            \ENDFOR
            \RETURN{$\left( a_{1}, a_{2}, \cdots, a_{n} \right)$}
        \end{algorithmic}
    \end{algorithm}
    \subsection{Static single-assignment form of bubble sort}
    \begin{algorithm}[H] 
        \caption{Static single-assignment form of bubble sort}
            \begin{algorithmic}[1]\label{algorithm_static_single_assignment_form_of_bubble_sort}
            \REQUIRE{$\left( a_{1}, a_{2}, \cdots, a_{n} \right)$}
            \ENSURE{a non-decreasing permutation of $\left( a_{1}, a_{2}, \cdots, a_{n} \right)$}
            \STATE{$\left( a_{1}^{\left( 0, n \right)}, a_{2}^{\left( 0, n \right)}, \cdots, a_{n}^{\left( 0, n \right)} \right) \gets \left( a_{1}, a_{2}, \cdots, a_{n} \right)$}
            \FOR{$i$ \textbf{from} $1$ \textbf{to} $n$}
                \STATE{$\left( a_{1}^{\left( i, 0 \right)}, a_{2}^{\left( i, 0 \right)}, \cdots, a_{n}^{\left( i, 0 \right)} \right) \gets \left( a_{1}^{\left( i - 1, n - i + 1 \right)}, a_{2}^{\left( i - 1, n - i + 1 \right)}, \cdots, a_{n}^{\left( i - 1, n - i + 1 \right)} \right)$}
                \FOR{$j$ \textbf{from} $1$ \textbf{to} $n - i$}
                    \IF{$a_{j}^{\left( i, j - 1 \right)} > a_{j + 1}^{\left( i, j - 1 \right)}$}
                        \STATE{$\left( a_{1}^{\left( i, j \right)}, a_{2}^{\left( i, j \right)}, \cdots, a_{n}^{\left( i, j \right)} \right) \gets \varsigma_{j} \left( a_{1}^{\left( i, j - 1 \right)}, a_{2}^{\left( i, j - 1 \right)}, \cdots, a_{n}^{\left( i, j - 1 \right)} \right)$}
                    \ELSE{}
                        \STATE{$\left( a_{1}^{\left( i, j \right)}, a_{2}^{\left( i, j \right)}, \cdots, a_{n}^{\left( i, j \right)} \right) \gets \left( a_{1}^{\left( i, j - 1 \right)}, a_{2}^{\left( i, j - 1 \right)}, \cdots, a_{n}^{\left( i, j - 1 \right)} \right)$}
                    \ENDIF
                \ENDFOR
            \ENDFOR
            \RETURN{$\left( a_{1}^{\left( n - 1, 1 \right)}, a_{2}^{\left( n - 1, 1 \right)}, \cdots, a_{n}^{\left( n - 1, 1 \right)} \right)$}
        \end{algorithmic}
    \end{algorithm}
    \subsection{Optimization of bubble sort}
    \begin{algorithm}[H]
        \caption{Optimization of bubble sort}\label{algorithm_optimization_of_bubble_sort}
        \begin{algorithmic}[1]
            \REQUIRE{$\left( a_{1}, a_{2}, \cdots, a_{n} \right)$}
            \ENSURE{a non-decreasing permutation of $\left( a_{1}, a_{2}, \cdots, a_{n} \right)$}
            \FOR{$i$ \textbf{from} $1$ \textbf{to} $n$}
                \STATE{$\xi \gets \mathcal{F}$}
                \FOR{$j$ \textbf{from} $1$ \textbf{to} $n - i$}
                    \IF{$a_{j} > a_{j + 1}$}
                        \STATE{swap $a_{j}, a_{j + 1}$}
                        \STATE{$\xi \gets \mathcal{T}$}
                    \ENDIF
                \ENDFOR
                \IF{$\xi = \mathcal{T}$}
                    \RETURN{$\left( a_{1}, a_{2}, \cdots, a_{n} \right)$}
                \ENDIF
            \ENDFOR
        \end{algorithmic}
    \end{algorithm}
    \subsection{Variant of optimization of bubble sort}
    \begin{algorithm}[H]
        \caption{Variant of optimization of bubble sort}
        \begin{algorithmic}[1]\label{algorithm_variant_of_optimization_of_bubble_sort}
            \REQUIRE{$\left( a_{1}, a_{2}, \cdots, a_{n} \right)$}
            \ENSURE{a non-decreasing permutation of $\left( a_{1}, a_{2}, \cdots, a_{n} \right)$}
            \FOR{$i$ \textbf{from} $1$ \textbf{to} $n$}
                \STATE{$\xi \gets \mathcal{F}$}
                \FOR{$j$ \textbf{from} $1$ \textbf{to} $n - i$}
                    \IF{$a_{j} > a_{j + 1}$}
                        \STATE{swap $a_{j}, a_{j + 1}$}
                        \STATE{$\xi \gets \mathcal{T}$}
                        \BREAK
                    \ENDIF
                \ENDFOR
                \IF{$\xi = \mathcal{T}$}
                    \RETURN{$\left( a_{1}, a_{2}, \cdots, a_{n} \right)$}
                \ENDIF
                \FOR{$k$ \textbf{from} $j + 1$ \textbf{to} $n - i$}
                    \IF{$a_{k} > a_{k + 1}$}
                        \STATE{swap $a_{k}, a_{k + 1}$}
                    \ENDIF
                \ENDFOR
            \ENDFOR
        \end{algorithmic}
    \end{algorithm}

\newpage
    
\section{Proofs}
    \subsection{Proof of theorem \ref{theorem_generalization_of_poisson}}
        \begin{proof}\label{proof_generalization_of_poisson_theorem}
            Using Cauchy-Schwarz inequality gives
            \begin{equation*}
                \left| \mathscr{S}_{n} \right| \sum_{\left\{ i, j \right\} \in \mathscr{S}_{n}} \left( \mathbb{E} \left( D_{n, \left\{ i, j \right\}} \right) \right)^{2} \geqslant \left( \sum_{\left\{ i, j \right\} \in \mathscr{S}_{n}} \mathbb{E} \left( D_{n, \left\{ i, j \right\}} \right) \right)^{2}
            \end{equation*}
            
            Note that $D_{n, \left\{ i, j \right\}} = D_{n, \left\{ i, j \right\}}^{2}$, then using $\mathbb{V} \left( W \right) = \mathbb{E} \left( W^{2} \right) - \left( \mathbb{E} \left( W \right) \right)^{2}$ gives
            \begin{equation*}
                \sum_{\left\{ i, j \right\} \in \mathscr{S}_{n}} \left( \mathbb{E} \left( D_{n, \left\{ i, j \right\}} \right) \right)^{2} = \mathbb{E} \left( \sum_{\left\{ i, j \right\} \in \mathscr{S}_{n}} D_{n, \left\{ i, j \right\}} \right) - \mathbb{V} \left( \sum_{\left\{ i, j \right\} \in \mathscr{S}_{n}} D_{n, \left\{ i, j \right\}} \right)
            \end{equation*}
            
            Therefore, $\displaystyle \mathbb{E} \left( \sum_{\left\{ i, j \right\} \in \mathscr{S}_{n}} D_{n, \left\{ i, j \right\}} \right) \to \lambda$ and $\displaystyle \mathbb{V} \left( \sum_{\left\{ i, j \right\} \in \mathscr{S}_{n}} D_{n, \left\{ i, j \right\}} \right) \to \lambda$ implies $\left| \mathscr{S}_{n} \right| \to \infty$ and $\left| T_{n} \right| \to \infty$.

            Based on Poisson limit theorem presented in \cite{Poisson_Approximation}, if $W_{n}$ follows $\mathrm{Poi} \left( \mu_{n} \right)$ where $\displaystyle \mu_{n} = \mathbb{E} \left( \sum_{\left\{ i, j \right\} \in \mathscr{S}_{n}} D_{n, \left\{ i, j \right\}} \right)$, then
            \begin{equation*}
                \begin{aligned}
                    & \quad \frac{1}{2} \sum_{k \in \mathbf{N}} \left| \mathbb{P} \left\{ \sum_{\left\{ i, j \right\} \in \mathscr{S}_{n}} D_{n, \left\{ i, j \right\}} = k \right\} - \mathbb{P} \left\{ W_{n} = k \right\} \right| \\
                    & \leqslant \frac{1 - \mathrm{e}^{- \mu_{n}}}{\mu_{n}} \sum_{\left\{ i, j \right\} \in \mathscr{S}_{n}} \left( \left( \mathbb{E} \left( D_{n, \left\{ i, j \right\}} \right) \right)^{2} + \sum_{\substack{\left\{ l, r \right\} \in \mathscr{S}_{n} \setminus \left\{ i, j \right\} \\ \left\{ l, r \right\} \cap \left\{ i, j \right\} \ne \varnothing}} \left( \mathbb{E} \left( D_{n, \left\{ i, j \right\}} \right) \mathbb{E} \left( D_{n, \left\{ l, r \right\}} \right) + \mathbb{E} \left( D_{n, \left\{ i, j \right\}} D_{n, \left\{ l, r \right\}} \right) \right) \right)
                \end{aligned}
            \end{equation*}

            Note that
            \begin{equation*}
                \begin{aligned}
                    & \quad \sum_{\left\{ i, j \right\} \in \mathscr{S}_{n}} \sum_{\substack{\left\{ l, r \right\} \in \mathscr{S}_{n} \setminus \left\{ i, j \right\} \\ \left\{ l, r \right\} \cap \left\{ i, j \right\} \ne \varnothing}} \mathbb{E} \left( D_{n, \left\{ i, j \right\}} \right) \mathbb{E} \left( D_{n, \left\{ l, r \right\}} \right) \\
                    & = \sum_{i \in \mathscr{T_{n}}} \left( \sum_{j \in T_{n} \setminus \left\{ i \right\}} \mathbb{E} \left( D_{n, \left\{ i, j \right\}} \right) \right)^{2} - 2 \sum_{\left\{ i, j \right\} \in \mathscr{S}_{n}} \left( \mathbb{E} \left( D_{n, \left\{ i, j \right\}} \right) \right)^{2} \\
                    & \leqslant \sum_{i \in \mathscr{T_{n}}} \left| T_{n} \setminus \left\{ i \right\} \right| \sum_{j \in T_{n} \setminus \left\{ i \right\}} \left( \mathbb{E} \left( D_{n, \left\{ i, j \right\}} \right) \right)^{2} - 2 \sum_{\left\{ i, j \right\} \in \mathscr{S}_{n}} \left( \mathbb{E} \left( D_{n, \left\{ i, j \right\}} \right) \right)^{2} \\
                    & = 2 \left( \left| T_{n} \right| - 2 \right) \sum_{\left\{ i, j \right\} \in \mathscr{S}_{n}} \left( \mathbb{E} \left( D_{n, \left\{ i, j \right\}} \right) \right)^{2}
                \end{aligned}
            \end{equation*}

            And there holds
            \begin{equation*}
                \sum_{\left\{ i, j \right\} \in \mathscr{S}_{n}} \sum_{\substack{\left\{ l, r \right\} \in \mathscr{S}_{n} \setminus \left\{ i, j \right\} \\ \left\{ l, r \right\} \cap \left\{ i, j \right\} \ne \varnothing}} \mathbb{E} \left( D_{n, \left\{ i, j \right\}} D_{n, \left\{ l, r \right\}} \right) =\sum_{i \in T_{n}} \sum_{j \in T_{n} \setminus \left\{ i \right\}} \sum_{k \in T_{n} \setminus \left\{ i, j \right\}} \mathbb{E} \left( D_{n, \left\{ i, j \right\}} D_{n, \left\{ i, k \right\}} \right)
            \end{equation*}

            Therefore, $\displaystyle \sum_{\left\{ i, j \right\} \in \mathscr{S}_{n}} D_{n, \left\{ i, j \right\}}$ converges in distribution to $\mathrm{Poi} \left( \lambda \right)$.
        \end{proof}
\end{document}